\documentclass[]{article}
\usepackage[dvips]{graphicx}
\usepackage{amsfonts}

\title{WORKSHOP\\
\vspace{1cm}
\bf SCIENCE and APPLICATIONS\\
of ADVANCED COMPUTING PARADIGMS\\
\vspace{1cm}
\it Centre of Excellence MIUR\\
Universit\`{a} di Padova\\
October, 28 -  29, 2004}
\author{}
\date{}

\begin{document}

\maketitle

\noindent

\section*{}

\begin{center}
{\tiny Padova, Aula Magna DEI, October 28, 2004}\\
\vspace{2cm}
{\bf \Large Computational Aspects of a Numerical Model for Combustion Flow}\\
\vspace{1cm}
{\bf {\large Gianluca Argentini}}\\
\vspace{0.2cm}
\normalsize gianluca.argentini@riellogroup.com \\
\vspace{0.2cm}
{\large Advanced Computing Laboratory}\\
\vspace{0.2cm}
{\it Riello Group}, Legnago (Verona), Italy
\end{center}

\newpage

\section*{Position of the problem}

Design, development and engineering of industrial power burners have 
strong mathematical requests:
\begin{itemize}
	\item numerical resolution of a PDEs system involving {\it Navier-Stokes} 
	equations for velocity and pressure fields, {\it energy conservation} law
	for temperature field, {\it Fick}'s law for diffusion of all the chemical
	species in the combustion chamber;
	\item geometrical design of the combustion head for a correct shape
	and optimal efficiency of {\it flame};
	\item geometrical design of {\it ventilation fans} and computation of 
	a correct air inflow for optimal combustion.
\end{itemize}

\begin{figure}[h]
	\begin{center}
	\includegraphics[width=6cm]{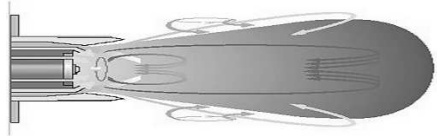}
	\caption{Combustion head and chamber for burner.}
	\end{center}
\end{figure}

\section*{Computational complexity analysis\\for~a~flow {\small (1)}}

{\bf {\it Simple example}} for a detailed knowledge of the velocity field
of fluid particles in the combustion chamber:

\begin{itemize}
	\item {\bf M} is the number of flow streamlines to compute;
	\item {\bf S} is the number of geometrical points for every streamline.
\end{itemize}
High values for {\bf M} are important for a {\it realistic simulation}
of the flow, high values for {\bf S} are important for a fine
{\it graphic resolution}: minimal values are of order {\it O}($10^3-10^4$).\\

\noindent Suppose to use a 3D grid 10 x 10 x 1000 cm (hence {\bf M} = 100, {\bf S} = 1000),
a medium value {\it $v_i$} = 50 cm/sec for every cartesian component
of velocity vector field, and a space resolution {\it h} = 0.5 cm.

\newpage

\section*{Computational complexity analysis\\for~a~flow {\small (2)}}

For numeric resolution of time-dependent advective PDEs, the {\it Courant-Friedrichs-Lewy}
({\it CFL}) {\it condition} gives an upper limit for the time step:\\

$\Delta t \leq \frac{ch}{v}$ \\

\noindent where {\it c} is a costant, usually $\leq 1$, depending on the
used numeric method, and ${\it v} = \sup |{\it v_j}|$.
The quantity $\frac{v \Delta t}{h}$ is called {\it CFL number}. 
Let {\it c} = 1; then\\

$\Delta t \leq \frac{0.5 cm}{50 \frac{cm}{sec}} = 0.01 \; sec$. \\

\noindent As consequence, for {\bf 1} real minute of simulation the flops are
of order {\bf {\it O}($10^{10}$)} and the occupation of RAM is {\bf {\it O}($10^0$) GB}:\\
{\it the computation is CPU expensive, RAM consuming and produces a lot of unuseful data}
(100 snapshots of the flow every second).

\section*{A Finite Differences method and Interpolations}
In the effort of minimize the relevance of these problems, we have studied 
a numeric model based on

\begin{itemize}
	\item a Finite Differences schema with a not too restrictive CFL condition;
	\item an appropriate interpolation of the numeric FD velocity-field for a finer
	resolution without modifying the grid step.
\end{itemize}

\noindent This model gives a numeric solution comparable with the solutions based
on finer grids: we present an {\it estimate} of its goodness and a mathematical 
justification.\\
The FD method is based on {\it Lax-Friedrichs} schema:

\begin{itemize}
	\item discretization in time: ${\partial_t u_j^n} = \frac{1}{\Delta t}(u_j^{n+1} - u_j^n)$,\\
	where $u_j^n \leftarrow \frac{1}{3}(u_{j+1}^n + u_j^n + u_{j-1}^n)$\\
	(for a better approximation we compute the mean on three values, 
	two in LF original form);
	\item discretization in space: ${\partial_x u_j^n}  = \frac{1}{2h}(u_{j+1}^n - u_{j-1}^n)$;
\end{itemize} 

\noindent where {\it u} is a velocity component, {\it n} the time step, {\it j} a value
on the cartesian coordinate {\it x}.

\newpage

\section*{Computational aspects of Lax-Friedrichs schema}
For this schema the CFL condition has costant {\it c} = 1; 
{\it the Finite Elements method with the same schema for discretization in time
has a more restrictive costant c} $< 1$.\\

\noindent If {\it K} $\in \mathbb{R^+}$, {\it K} $\leq \frac{1}{2}$, we can
define the norm ~$\parallel${\it u}$\parallel$ = {\it K} $sup_j |{\it u_j}|$;
then the modified LF schema is {\it strongly stable}:
$\parallel$${\it u^n}$$\parallel$ $\leq$ $\parallel$${\it u^{n-1}}$$\parallel$ $\forall {\it n} \in \mathbb{N}$;
hence there is not the {\it blowing up} of the numeric solution.

\noindent Suppose we want to compute at most 10 snapshots for every second;
then, in the hypothesis {\it v} = 50 cm/sec as the previous example, from\\

$v \Delta t \leq h$\\

\noindent we must use as minimum a grid step {\it h} = 5 cm.\\
\noindent This case gives {\bf S} = 200, the total flops for 1 minute
of simulation is now of order {\it O}($10^8$) and the occupation of RAM
is of order {\it O}($10^{-2}$) GB.\\
The gain is of order {\it O}($10^2$).\\
The grid step {\it h} = 5 cm is too big for a good resolution of
streamlines for flows into the combustion head: for better final results,
it can be useful a method based on {\it interpolations} of the computed LF values.

\section*{Interpolation of trajectories {\small (1)}}
Every streamline of LF solution is divided into {\bf N} couples of points,\\
\{$(P_1,P_2),(P_2,P_3),...,(P_{{\bf N}-1},P_{\bf N})$\}, so that {\bf S} = {\bf N}+1.\\
We use for every couple a cubic polynomial ({\it spline}) imposing
the following four analytical conditions ({\bf v} is the LF solution):

\begin{itemize}
	\item passage at ${\it P_k}$ point, $1 \leq {\it k} \leq {\bf N}-1$;
	\item passage at ${\it P_{k+1}}$ point;
	\item the first derivative at ${\it P_k}$ is equal to ${\bf v_{\it k}}$;
	\item the first derivative at ${\it P_{k+1}}$ is equal to ${\bf v_{\it {k+1}}}$.
\end{itemize}

\noindent In this way we can construct a set of class ${\it C}^1$ 
new trajectories; we want to estimate

\begin{enumerate}
	\item the overload for finding and valuating all the cubics;
	\item the difference compared to the real LF solution of the smaller grid step.
\end{enumerate}

\newpage

\section*{Interpolation of trajectories {\small (2)}}
For simplicity, consider a single component of a cubic:\\
$s(t) = at^3 + bt^2 + ct + d$, where $0 \leq t \leq 1$;\\
if {\bf T} is the $4\times4$ matrix\\

${\bf T} = \left( \begin{array}{cccc}
2  & -2 & 1  & 1 \\
-3 & 3  & -2 & -1 \\
0  & 0  & 1  & 0 \\
1  & 0  & 0  & 0 
\end{array} \right)$\\

\noindent and $(p_1, p_2, v_1, v_2)$ is the vector of cartesian coordinates and
velocities components of points $P_1$ and $P_2$, we have\\

$(a, b, c, d) = {\bf T} (p_1, p_2, v_1, v_2)$.

\section*{Interpolation of trajectories {\small (3)}}
We define the $4{\bf M}\times4{\bf M}$ {\it global matrix}\\

${\bf G} = \left( \begin{array}{ccccc}
{\bf T}  & {\bf 0} & . & . &  {\bf 0} \\
{\bf 0}  & {\bf T} & . & . &  {\bf 0} \\
. & . & . & . & . \\
. & . & . & . & . \\
{\bf 0}  & {\bf 0} & . & . &  {\bf T}
\end{array} \right)$\\

\noindent where {\bf 0} is the $4\times4$ zero-matrix. Then

\begin{itemize}
	\item {\bf G} is a {\it sparse} matrix with density number $\leq \frac{1}{\bf M}$;
	\item if ${\bf p} = (p_{(1,1)}, p_{(1,2)}, . . ., v_{({\bf M},1)}, v_{({\bf M},2)})$,
	we can compute the cubics, between two points, for all the {\bf M} trajectories 
	by the product {\bf G}{\bf p}.
\end{itemize}

\newpage

\section*{Interpolation of trajectories {\small (4)}}
The theoric number of flops for computing the coefficients of all the splines is
of order ${\it O}(10 {\bf M}^2{\bf N})$. If {\bf M} = $10^4$ and {\bf N} = $10^3$,
the total number of flops is ${\it O}(10^{12})$.\\

\noindent With a single processor having a clock frequency of {\it O}(1) {\bf GHz},
the total time can require some hundreds of seconds, a performance not very good
for practical purposes; using

\begin{itemize}
	\item some mathematical libraries as {\it LAPACK} routines with {\bf Fortran}
	calls or {\bf Matlab} environment,
	\item distributed computation on a multinode cluster,
\end{itemize}

\noindent we have reached a computation time of some tens of seconds.\\

\noindent {\it Example}: {\bf Matlab} has internal Lapack level 3 {\bf BLAS} routines
for fast matrix-matrix multiplication and treatment of sparse matrices.

\section*{Interpolation of trajectories {\small (5)}}

\begin{figure}[h]
	\begin{center}
	\includegraphics[width=6.5cm]{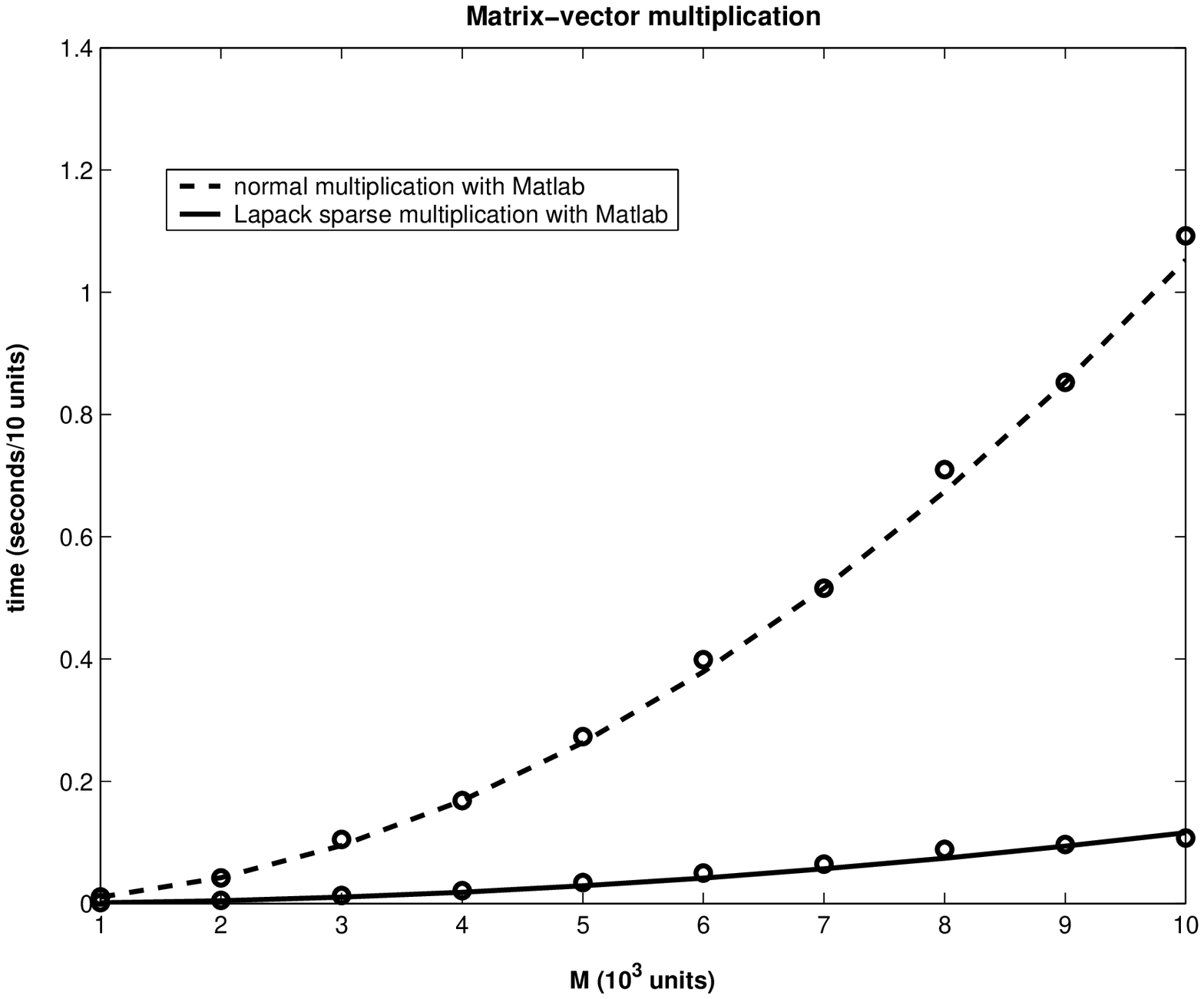}
	\end{center}
\end{figure}

\noindent Performances for a single {\bf Gp} multiplication using an Intel
Xeon 3.2 GHz with 1 MB internal {\it cache}: for {\bf M}=$10^4$ the memory occupied
by the sparse version of {\bf G} is only {\it O}($10^2$) KB instead of theoric
{\it O}($10^6$): {\bf G} can be stored in processor cache.

\newpage

\section*{Computation of splines values {\small (1)}}
Now we need a fast method for computing the splines values in a set of
{\it parameter~ticks} with fine sampling.\\
Let ${\it r} \in \mathbb{N}^+$ the number of ticks for each cubic: then the values of the
parameter {\it t} in these ticks are $(0, \frac{1}{r}, \frac{2}{r}, . . ., \frac{r-1}{r}, 1)$;
the value of a cubic at ${\it t_0}$ is a {\it scalar product}:\\

${at_0^3 + bt_0^2 + ct + d = (a, b, c, d)\cdot(t_0^3, t_0^2, t_0, 1)}$.\\

\noindent Consider the constant $4\times{\it (r+1)}$ {\bf R} matrix
and the (${\bf M}\times4$) {\bf C} matrix:\\

${\bf R} = \left( \begin{array}{cccccc}
0 & (\frac{1}{r})^3 & . & . & (\frac{r-1}{r})^3 & 1 \\
0 & (\frac{1}{r})^2 & . & . & (\frac{r-1}{r})^2 & 1 \\
0 & (\frac{1}{r})^1 & . & . & (\frac{r-1}{r})^1 & 1\\
1 & 1 & . & . & 1 & 1
\end{array} \right)$\\

${\bf C} = \left( \begin{array}{cccc}
a_1 & b_1 & c_1 & d_1 \\
a_2 & b_2 & c_2 & d_2 \\
. & . & . & . \\
. & . & . & . \\
a_{\bf M} & b_{\bf M} & c_{\bf M} & d_{\bf M}
\end{array} \right)$\\

\section*{Computation of splines values {\small (2)}}
Then the ${\bf M}\times{\it (r+1)}$ matrix {\bf C}{\bf R} contains for each row the
values of a cubic between two points, for all the trajectories ({\it eulerian}
method: computation of all the position and velocity at a fixed instant). The
flops for one multiplication are of order {\it O}(10{\bf M}{\it r}).\\

\noindent Tests with Xeon 3.2 GHz processor, {\bf M} = $10^4$, {\it r} = 10 and
{\it GNU} {\bf Fortran77} show a time of {\it 0.01 seconds} for a multiplication.\\
With {\bf N} = $10^2$, the time for computing the values of all the splines
of a single time step is {\it 4.5} seconds (theoric for 3D: $0.01\times10^2\times3 = {\it 3}$ secs).\\

\noindent If {\bf p} is the number of available processors and {\it mod}({\bf M}, {\bf p}) = 0,
the computation can be parallelized distributing $\frac{\bf M}{\bf p}$ rows
of matrix {\bf C} to each processor: there is no need of communication among
processes.\\
A version of {\bf High Performance Fortran} on a SMP system with 4 ItaniumII
processors shows a quasilinear speedup for {\bf M}, {\bf N} of order {\it O}($10^3$).

\newpage

\section*{Time for computation}

\begin{figure}[h]
	\begin{center}
	\includegraphics[width=7cm]{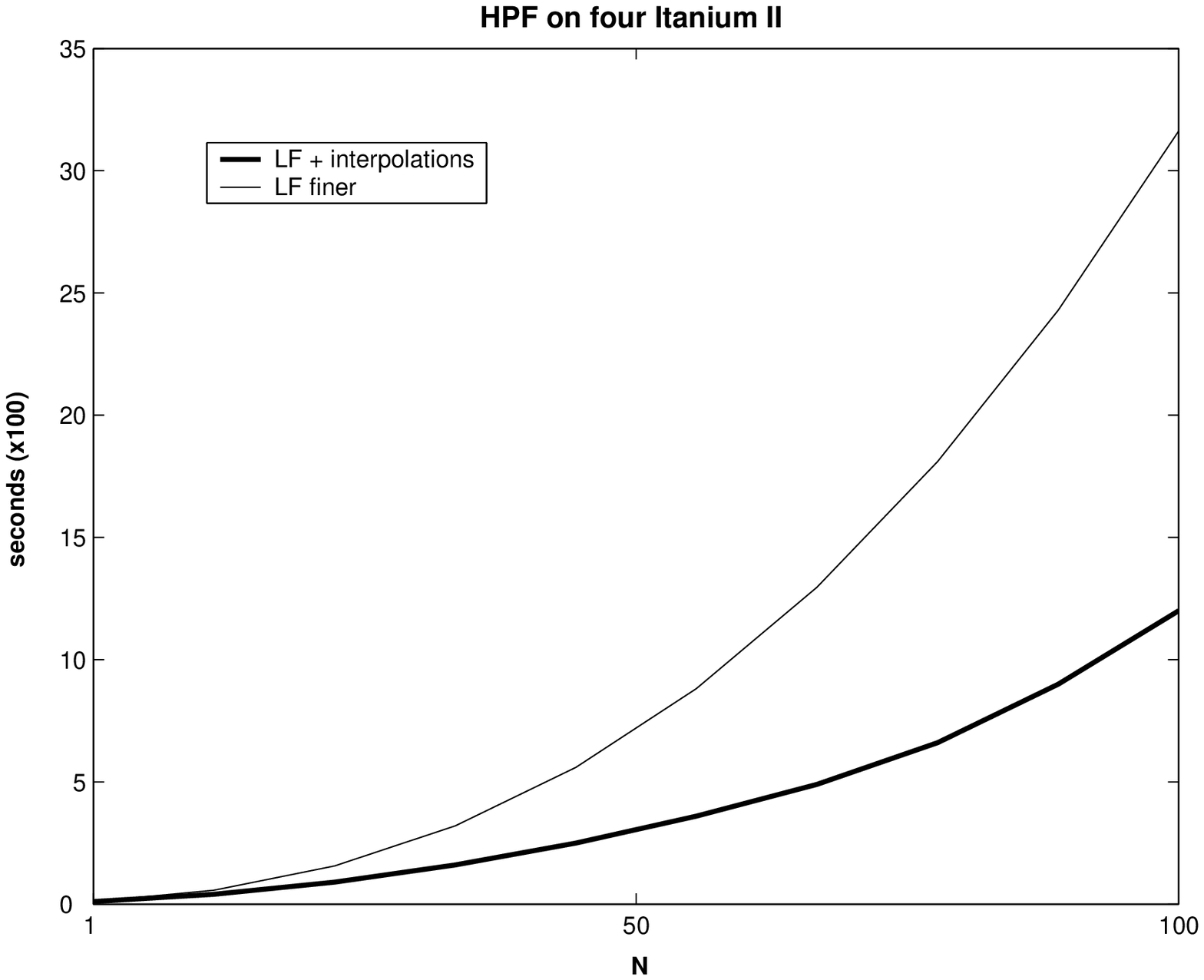}
	\end{center}
\end{figure}

\noindent These are the total time of computation for the two methods in the case of a cilinder of length
{\bf L} = 1m, a flow with a max. speed {\it v} = 10cm/sec, {\bf M} = $10^4$, {\it r} = 10 and 1 minute
of real simulation. The space grid is {\it h} = $\frac{{\bf L}}{10{\bf N}}$.

\section*{Estimate of LF+interpolations vs normal LF {\small (1)}}
{\it But what is the difference between the modified LF solution and
normal LF solution?}\\
Consider the one-dimensional case.
Let {\bf u}=$(u_k)$ the solution of normal LF schema with grid step {\it h} and
initial value {$\bf u_0$}; {\bf w}=$(w_m)$ the solution of normal LF schema 
with grid step {\it $s\times h$}, $s \in \mathbb{N}^+$, and initial value {$\bf w_0$} $\subset$ {$\bf u_0$};
{\bf v}=$(v_n)$ the solution of modified LF schema obtained by interpolation
of {\bf w} and valuation on {\it s} points per cubic; for a cubic, let ${v_k}$,
$k \leq s$, the value of {\bf v} at {\it t}=$\frac{k}{s}$ and ${u_k}$ the
value of {\it u} at the corresponding node of the finer grid; $\frac{v \Delta t}{h}$ the
CFL number and {\it N} the {\it N}-th time step. Let
$${\it M_0} = \max_{|m-n|=1}|u_{0,m}-u_{0,n}|$$
Then it is possible to prove this result:
\newtheorem{teorema}{Theorem}
\begin{teorema}
If ${\it M_0} >$ {\normalsize 0}, there are two positive constants {\it A} and {\it B} such that\\
$$|v_n - u_n| \leq (A + Bs)M_0\sum_{i=0}^{\it N}(\frac{v \Delta t}{2h})^i\qquad \forall n\in\{grid\;indexes\}, \forall N\in\mathbb{N}.$$
\end{teorema}

\newpage

\section*{Estimate of LF+interpolations vs normal LF {\small (2)}}
The CFL number $\frac{v \Delta t}{h}$ is usually indicated by $\lambda v$.
From the previous theorem it follows:
\newtheorem{corollario}{Corollary}
\begin{corollario}
If $\lambda v < 2$, then\\
$$|v_n - u_n| < \frac{2(A + Bs)M_0}{2 - \lambda v}.$$
\end{corollario}
\noindent The CFL condition satisfies the hypothesis of the corollary.\\
Hence, for a realistic solution from the LF+interpolations model,
the conditions are:

\begin{itemize}
	\item a small $(\ll 2)$ CFL number,
	\item a not too big number {\it s} of valuations for the cubics; note that 
	{\it s} has the inverse logical meaning of the previous {\it r} parameter.
\end{itemize}

\noindent Note that if ${\it M_0}$ is very big, as in the case of very caotic
flows, the LF+interpolations solution can be not very realistic.

\section*{Estimate of LF+interpolations vs normal LF {\small (3)}}
Testing the estimate: {\it example} for one-dimensional non linear Navier-Stokes equation,
$\lambda v$ = 1, {\it s} = 10, after {\it N} = $10^5$ time steps; graphic of the error between
LF+interpolations and normal LF solutions.

\begin{figure}[h]
	\begin{center}
	\includegraphics[width=6.0cm]{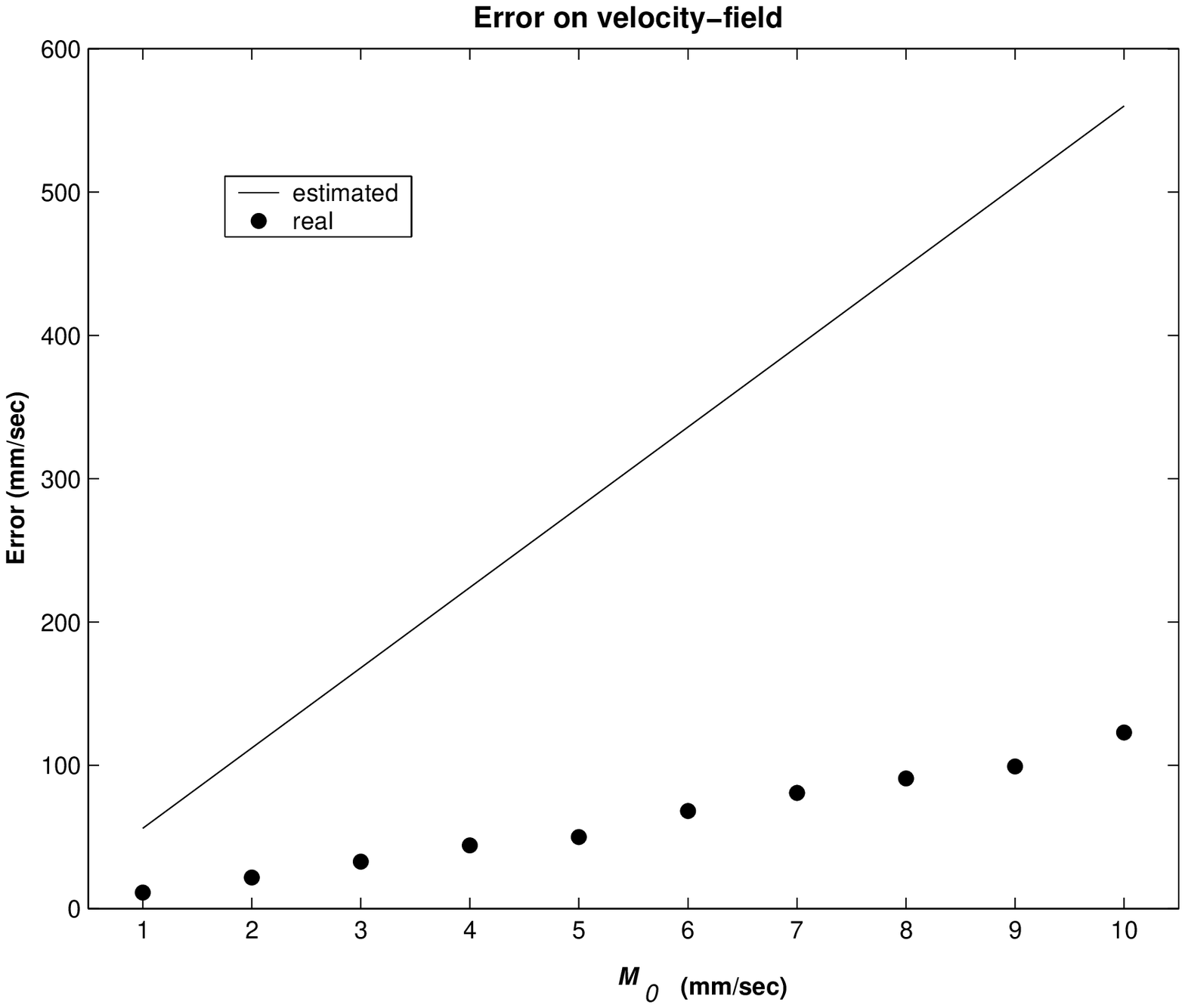}
	\end{center}
\end{figure}

\noindent In this case it can be shown that {\it A}~= 8, {\it B} = 2 is a 
first, not optimized, approximation for the two constants.
The picture shows that the estimate is correct but large.

\newpage

\section*{Conclusions}
The numeric LF schema can be modified using the interpolations
method so that:

\begin{itemize}
	\item the time spent on computation is much lower than the time of the LF based
on the corresponding finer grid;
	\item the computation can be parallelized on multiprocessors environment with
	very reduced need of communication;
	\item the error on normal LF solution can be estimated and depends on the 
	initial value ${\bf u_0}$ of the problem;
	\item the estimate is compatible with CFL condition.
\end{itemize}

$\;$ \\
$\;$ \\
\noindent {\textbf{\textit{Thanks}}

\end{document}